\documentclass[11pt]{amsart}
\usepackage[margin=30mm]{geometry}
\usepackage{amsmath,amssymb}
\usepackage{amsthm}
\usepackage{mathrsfs}

\newtheorem{thm}{Theorem}[section]

\newtheorem{lem}{Lemma}[section]
\newtheorem{cor}{Corollary}[section]
\newtheorem{defi}{Definition}[section]

\newtheorem{rem}{Remark}[section]

\usepackage{enumitem}

\begin{document}

\title{Shadowing and the basins of terminal chain components}
\author{Noriaki Kawaguchi}
\subjclass[2020]{37B65}
\keywords{shadowing; basin; chain component; generic; chain continuous}
\address{Research Institute of Science and Technology, Tokai University, 4-1-1 Kitakaname, Hiratsuka, Kanagawa 259-1292, Japan}
\email{gknoriaki@gmail.com}

\begin{abstract}
We provide an alternative view of some results in \cite{A1,AHK,H}. In particular, we prove that (1) if a continuous self-map of a compact metric space has the shadowing, then the union of the basins of terminal chain components is a dense $G_\delta$-subset of the space; and (2) if a continuous self-map of a locally connected compact metric space has the shadowing, and if the chain recurrent set is totally disconnected, then the map is almost chain continuous.   
\end{abstract}

\maketitle

\markboth{NORIAKI KAWAGUCHI}{Shadowing and the basins of terminal chain components}

\section{Introduction}

{\em Shadowing} is an important concept in the topological theory of dynamical systems (see \cite{AH,P2} for background). It was derived from the study of hyperbolic differentiable dynamics \cite{A,B} and generally refers to a situation in which coarse orbits, or {\em pseudo-orbits}, can be approximated by true orbits. Above all else, it is worth mentioning that the shadowing is known to be {\em generic} in the space of homeomorphisms or continuous self-maps of a closed differentiable manifold \cite{MO,PP} and so plays a significant role in the study of topologically generic dynamics.

{\em Chain components} are basic objects for global understanding of dynamical systems \cite{C}.  In this paper, we focus on attractor-like, or  {\em terminal}, chain components and the basins of them. By a result (Corollary 6.16) of \cite{H}, if a continuous flow on a compact metric space has the so-called {\em weak shadowing}, then the union of the basins of terminal chain components is a dense $G_\delta$-subset of the space. For any continuous self-map of a compact metric space, we strengthen it by assuming the standard shadowing (Theorem 1.1). Our proof is by a method related to but independent of a result (Proposition 22 in Section 7) of \cite{A1}. It is shown in \cite{AHK} that topologically generic homeomorphisms of a closed differentiable manifold are almost chain continuous. We also give an alternative proof of this fact by using the genericity of shadowing.

First we define the chain components. Throughout, $X$ denotes a compact metric space endowed with a metric $d$. 

\begin{defi}
\normalfont
Given a continuous map $f\colon X\to X$ and $\delta>0$, a finite sequence $(x_i)_{i=0}^{k}$ of points in $X$, where $k>0$ is a positive integer, is called a {\em $\delta$-chain} of $f$ if $d(f(x_i),x_{i+1})\le\delta$ for every $0\le i\le k-1$.  A $\delta$-chain $(x_i)_{i=0}^{k}$ of $f$ with $x_0=x_k$ is said to be a {\em $\delta$-cycle} of $f$. 
\end{defi}

Let $f\colon X\to X$ be a continuous map. For any $x,y\in X$ and $\delta>0$, the notation $x\rightarrow_\delta y$ means that there is a $\delta$-chain $(x_i)_{i=0}^k$ of $f$ with $x_0=x$ and $x_k=y$. We write $x\rightarrow y$ if $x\rightarrow_\delta y$ for all $\delta>0$. We say that $x\in X$ is a {\em chain recurrent point} for $f$ if $x\rightarrow x$, or equivalently, for every $\delta>0$, there is a $\delta$-cycle $(x_i)_{i=0}^{k}$ of $f$ with $x_0=x_k=x$. Let $CR(f)$ denote the set of chain recurrent points for $f$. We define a relation $\leftrightarrow$ in
\[
CR(f)^2=CR(f)\times CR(f)
\]
by: for any $x,y\in CR(f)$, $x\leftrightarrow y$ if and only if $x\rightarrow y$ and $y\rightarrow x$. Note that $\leftrightarrow$ is a closed equivalence relation in $CR(f)^2$ and satisfies $x\leftrightarrow f(x)$ for all $x\in CR(f)$. An equivalence class $C$ of $\leftrightarrow$ is called a {\em chain component} for $f$.  We regard the quotient space
\[
\mathcal{C}(f)=CR(f)/{\leftrightarrow}
\]
as a space of chain components.

A subset $S$ of $X$ is said to be $f$-invariant if $f(S)\subset S$. For an $f$-invariant subset $S$ of $X$, we say that $f|_S\colon S\to S$ is {\em chain transitive} if for any $x,y\in S$ and $\delta>0$, there is a $\delta$-chain $(x_i)_{i=0}^k$ of $f|_S$ with $x_0=x$ and $x_k=y$.

\begin{rem}
\normalfont
The following properties hold:
\begin{itemize}
\item $CR(f)=\bigsqcup_{C\in\mathcal{C}(f)}C$,
\item every $C\in\mathcal{C}(f)$ is a closed $f$-invariant subset of $CR(f)$,
\item $f|_C\colon C\to C$ is chain transitive for all $C\in\mathcal{C}(f)$,
\item for any $f$-invariant subset $S$ of $X$, if $f|_S\colon S\to S$ is chain transitive, then $S\subset C$ for some $C\in\mathcal{C}(f)$.
\end{itemize}
\end{rem}

Next we recall the definition of terminal chain components.

\begin{defi}
\normalfont
We say that a closed $f$-invariant subset $S$ of $X$ is {\em chain stable} if for any $\epsilon>0$, there is $\delta>0$ such that every $\delta$-chain $(x_i)_{i=0}^k$ of $f$ with $x_0\in S$ satisfies $d(x_i,S)\le\epsilon$ for all $0\le i\le k$. Following \cite{AHK}, we say that $C\in\mathcal{C}(f)$ is {\em terminal} if $C$ is chain stable. We denote by $\mathcal{C}_{\rm ter}(f)$ the set of terminal chain components for $f$.
\end{defi}

\begin{rem}
\normalfont
For any continuous map $f\colon X\to X$, a partial order $\le$ on $\mathcal{C}(f)$ is defined by: for all $C,D\in\mathcal{C}(f)$, $C\le D$ if and only if $x\rightarrow y$ for some $x\in C$ and $y\in D$. We can easily show that for any $C\in\mathcal{C}(f)$, $C\in\mathcal{C}_{\rm ter}(f)$ if and only if $C$ is maximal with respect to $\le$, that is, $C\le D$ implies $C=D$ for all $D\in\mathcal{C}(f)$. 
\end{rem}

Given a continuous map $f\colon X\to X$ and $x\in X$, the {\em $\omega$-limit set} $\omega(x,f)$ of $x$ for $f$ is defined as the set of $y\in X$ such that
\[
\lim_{j\to\infty}f^{i_j}(x)=y
\]
for some sequence $0\le i_1<i_2<\cdots$. Note that $\omega(x,f)$ is a closed $f$-invariant subset of $X$ and $f|_{\omega(x,f)}\colon\omega(x,f)\to\omega(x,f)$ is chain transitive. We denote by $C(x,f)$ the unique $C(x,f)\in\mathcal{C}(f)$ such that $\omega(x,f)\subset C(x,f)$. For each $C\in\mathcal{C}(f)$, we define the {\em basin} $W^s(C)$ of $C$ by
\[
W^s(C)=\{x\in X\colon\lim_{i\to\infty}d(f^i(x),C)=0\}.
\]
For every $x\in X$, since
\[
\lim_{i\to\infty}d(f^i(x),\omega(x,f))=0,
\]
we have $x\in W^s(C)$ if and only if $C=C(x,f)$. This implies
\[
\{x\in X\colon C(x,f)\in\mathcal{C}_{\rm ter}(f)\}=\bigsqcup_{C\in\mathcal{C}_{\rm ter}(f)}W^s(C).
\]
We also define the {\em chain $\omega$-limit set} $\omega^\ast(x,f)$ of $x$ for $f$ as the set of $y\in X$ such that for any $\delta>0$ and $N>0$, there is a $\delta$-chain $(x_i)_{i=0}^k$ of $f$ with $x_0=x$, $x_k=y$, and $k\ge N$. Note that $\omega^\ast(x,f)$ is a closed $f$-invariant subset of $X$ and chain stable. We have
\[
\omega(x,f)\subset C(x,f)\subset\omega^\ast(x,f).
\]

\begin{rem}
\normalfont
The chain $\omega$-limit set is denoted in \cite{AHK} as $\omega\mathcal{C}(x,f)$ instead of $\omega^\ast(x,f)$.
\end{rem}

The following lemma is obvious (see Section 1.4 of \cite{AHK}).

\begin{lem}
Let $f\colon X\to X$ be a continuous map.
\begin{itemize}
\item[(A)] For any $x\in X$, the following properties are equivalent:
\begin{itemize}
\item $C(x,f)\in\mathcal{C}_{\rm ter}(f)$,
\item $\omega^\ast(x,f)\subset C(x,f)$,
\item $\omega^\ast(x,f)=C(x,f)$,
\item $f|_{\omega^\ast(x,f)}\colon\omega^\ast(x,f)\to\omega^\ast(x,f)$ is chain transitive.
\end{itemize}
\item[(B)] For any $x\in X$, the following properties are equivalent:
\begin{itemize}
\item $\omega(x,f)=C(x,f)=\omega^\ast(x,f)$,
\item $C(x,f)\in\mathcal{C}_{\rm ter}(f)$ and $\omega(x,f)=C(x,f)$.
\end{itemize}
\end{itemize}
\end{lem}

We give the definition of shadowing.

\begin{defi}
\normalfont
Let $f\colon X\to X$ be a continuous map and let $\xi=(x_i)_{i\ge0}$ be a sequence of points in $X$. For $\delta>0$, $\xi$ is called a {\em $\delta$-pseudo orbit} of $f$ if $d(f(x_i),x_{i+1})\le\delta$ for all $i\ge0$. For $\epsilon>0$, $\xi$ is said to be {\em $\epsilon$-shadowed} by $x\in X$ if $d(f^i(x),x_i)\leq \epsilon$ for all $i\ge 0$. We say that $f$ has the {\em shadowing property} if for any $\epsilon>0$, there is $\delta>0$ such that every $\delta$-pseudo orbit of $f$ is $\epsilon$-shadowed by some point of $X$. 
\end{defi}

For a topological space $Z$, a subset $S$ of $Z$ is called a {\em $G_\delta$-subset} of $Z$ if $S$ is a countable intersection of open subsets of $Z$. If $Z$ is completely metrizable, then by Baire Category Theorem, every countable intersection of open dense subsets of $Z$ is dense in $Z$. We know that a subspace $Y$ of a completely metrizable space $Z$ is completely metrizable if only if $Y$ is a $G_\delta$-subset of $Z$ (see Theorem 24.12 of \cite{W}).

For any continuous map $f\colon X\to X$ and $x\in X$, let $\Omega(x,f)$ denote the set of $y\in X$ such that
\[
\lim_{j\to\infty}f^{i_j}(x_j)=y
\]
for some sequence $0\le i_1<i_2<\cdots$ and $x_j\in X$, $j\ge0$, with
\[
\lim_{j\to\infty}x_j=x.
\]
Note that
\[
\omega(x,f)\subset\Omega(x,f)\subset\omega^\ast(x,f)
\]
for all $x\in X$. By Proposition 22 in Section 7 of \cite{A1}, we know that
\[
\{x\in X\colon\omega(x,f)=\Omega(x,f)\}
\]
is a dense $G_{\delta}$-subset of $X$. The proof of this result in \cite{A1} is based on the fact the set of continuity points of a lower semicontinuous (lsc) set-valued map is a dense $G_{\delta}$-subset and it is not trivial. If $f$ has the shadowing property, we have 
\[
\Omega(x,f)=\omega^\ast(x,f)
\]
for all $x\in X$; therefore,
\[
\{x\in X\colon\omega(x,f)=C(x,f)=\omega^\ast(x,f)\}=\{x\in X\colon\text{$C(x,f)\in\mathcal{C}_{\rm ter}(f)$ and $\omega(x,f)=C(x,f)$}\}
\]
is a dense $G_{\delta}$-subset of $X$ (see \cite{H} and \cite{P1} for related results). The main aim of this paper is to give an alternative proof of the following statement.

\begin{thm}
If a continuous map $f\colon X\to X$ has the shadowing property, then
\[
V(f)=\{x\in X\colon C(x,f)\in\mathcal{C}_{\rm ter}(f)\} 
\]
and
\[
W(f)=\{x\in V(f)\colon\omega(x,f)=C(x,f)\} 
\]
are dense $G_\delta$-subsets of $X$.
\end{thm}

Given a continuous map $f\colon X\to X$ and $x\in X$, we say that $f$ is {\em chain continuous} at $x$ if for any $\epsilon>0$, there is $\delta>0$ such that every $\delta$-pseudo orbit $(x_i)_{i\ge0}$ of $f$ with $x_0=x$ is $\epsilon$-shadowed by $x$ \cite{A2}. We denote by $CC(f)$ the set of chain continuity points for $f$. We say that a closed $f$-invariant subset $S$ of $X$ is an {\em odometer} if $(S,f|_S)$ is topologically conjugate to an odometer. This is equivalent to that $S$ is a Cantor space and
\[
f|_{S}\colon S\to S
\]
is a minimal equicontinuous homeomorphism. By Theorem 7.5 of \cite{AHK}, we know that for any $x\in X$, $x\in CC(f)$ if and only if
\[
\omega(x,f)=C(x,f)=\omega^\ast(x,f)
\]
and $C(x,f)$ is a periodic orbit or an odometer. By Lemma 1.1, this is equivalent to that $C(x,f)\in\mathcal{C}_{\rm ter}(f)$ and $C(x,f)$ is a periodic orbit or an odometer. We say that $X$ is {\em locally connected} if for any $x\in X$ and any open subset $U$ of $X$ with $x\in U$, we have $x\in V\subset U$ for some open connected subset $V$ of $X$. If $X$ is locally connected and $CR(f)$ is totally disconnected, then due to results of \cite{BS} or \cite{HH}, every $C\in\mathcal{C}_{\rm ter}(f)$ is a periodic orbit or an odometer. By these facts, we obtain the following lemma.

\begin{lem}
Let $f\colon X\to X$ be a continuous map. If $X$ is locally connected and $CR(f)$ is totally disconnected, then for any $x\in X$, the following properties are equivalent:
\begin{itemize}
\item $x\in CC(f)$,
\item $\omega(x,f)=C(x,f)=\omega^\ast(x,f)$,
\item $C(x,f)\in\mathcal{C}_{\rm ter}(f)$.
\end{itemize} 
\end{lem}

Let $f\colon X\to X$ be a continuous map. For any $j,l\ge1$, let $C_{j,l}$ denote the set of $x\in X$ such that there is a neighborhood $U$ of $x$ for which every $\frac{1}{j}$-pseudo orbit $(x_i)_{i\ge0}$ of $f$ with $x_0\in U$ is $\frac{1}{l}$-shadowed by $x_0$. We see that $C_{j,l}$ is an open subset of $X$ for all $j,l\ge1$ and
\[
CC(f)=\bigcap_{l\ge1}\bigcup_{j\ge1}C_{j,l},
\]
thus $CC(f)$ is a $G_{\delta}$-subset of $X$. We say that $f$ is {\em almost chain continuous} if $CC(f)$ is a dense $G_\delta$-subset of $X$. By Theorem 1.1 and Lemma 1.2, we obtain the following theorem.

\begin{thm}
Let $f\colon X\to X$ be a continuous map. If $X$ is locally connected, $f$ has the shadowing property, and if $CR(f)$ is totally disconnected, then $f$ is almost chain continuous.
\end{thm}

We present a corollary of Theorem 1.2. For a closed differentiable manifold $M$, let $\mathcal{H}(M)$ (resp.\:$\mathcal{C}(M)$) denote the set of homeomorphisms (resp.\:continuous self-maps) of $M$, endowed with the $C^0$-topology. It is shown in \cite{AHK} that generic $f\in\mathcal{H}(M)$ (resp.\:$f\in\mathcal{C}(M)$, if $\dim{M}>1$) is almost chain continuous. Note that the shadowing is generic in $\mathcal{H}(M)$ \cite{PP} and also generic in $\mathcal{C}(M)$ \cite{MO}. Moreover, by results of \cite{AHK, KOU}, we know that for generic $f\in\mathcal{H}(M)$ (resp.\:$f\in\mathcal{C}(M)$), $CR(f)$ is totally disconnected. Thus, by Theorem 1.2, we obtain the following corollary.

\begin{cor}
Generic $f\in\mathcal{H}(M)$ (resp.\:$f\in\mathcal{C}(M)$) is almost chain continuous. 
\end{cor}

Our results also apply to the case where $X$ is not a manifold. We say that $X$ is a {\em dendrite} if $X$ is connected, locally connected, and contains no simple closed curves. The shadowing is proved to be generic in the space of continuous self-maps of a dendrite \cite{BMR,KMOK}. On the other hand, by results of \cite{KOU}, a generic continuous self-map of a dendrite has totally disconnected chain recurrent set. By Theorem 1.2, we conclude that a generic continuous self-map of a dendrite is almost chain continuous.

This paper consists of two sections. In the next section, we prove Theorem 1.1. 

\section{Proof of Theorem 1.1}

In this section, we prove Theorem 1.1. The proof is based on the following lemma in \cite{K}.

\begin{lem}[{\cite[Lemma 2.1]{K}}]
For any continuous map $f\colon X\to X$ and $x\in X$, there is $C\in\mathcal{C}_{\rm ter}(f)$ such that for every $\delta>0$, there is a $\delta$-chain $(x_i)_{i=0}^k$ of $f$ with $x_0=x$ and $x_k\in C$.
\end{lem}

By using this lemma, we prove Theorem 1.1.

\begin{proof}[Proof of Theorem 1.1]
First, we show that $V(f)$ is a dense $G_\delta$-subset of $X$. For any subset $A$ of $X$ and $r>0$, we denote by $U_r(A)$ the open $r$-neighborhood of $A$:
\[
U_r(A)=\{x\in X\colon d(x,A)<r\}.
\]
Fix a sequence $(\epsilon_j)_{j\ge1}$ of positive numbers such that $\epsilon_1>\epsilon_2>\cdots$ and
\[
\lim_{j\to\infty}\epsilon_j=0.
\]
For any $j\ge1$ and $C\in\mathcal{C}_{\rm ter}(f)$, we take $\delta_{j,C}>0$ such that $x\in U_{\delta_{j,C}}(C)$ implies
\[
\omega^\ast(x,f)\subset U_{\epsilon_j}(C)
\]
for all $x\in X$. Let
\[
U_{j,C}=U_{\delta_{j,C}}(C)
\]
for all $j\ge1$ and $C\in\mathcal{C}_{\rm ter}(f)$. We define a subset $V$ of $X$ by
\[
V=\bigcap_{j\ge1}\bigcup_{C\in\mathcal{C}_{\rm ter}(f)}\bigcup_{m\ge0}f^{-m}(U_{j,C}).
\]
Note that $V$ is a $G_{\delta}$-subset of $X$. Since $f$ has the shadowing property, by Lemma 2.1, we see that for every $x\in X$, there is $C\in\mathcal{C}_{\rm ter}(f)$ such that
\[
x\in\overline{\bigcup_{m\ge0}f^{-m}(U_{j,C})}
\]
for all $j\ge1$. This implies that $V$ is a dense $G_\delta$-subset of $X$. It remains to prove that $V(f)=V$. Given any $x\in V(f)$, by $C(x,f)\in\mathcal{C}_{\rm ter}(f)$ and
\[
x\in\bigcap_{j\ge1}\bigcup_{m\ge0}f^{-m}(U_{j,C(x,f)})\subset V,
\]
we have $x\in V$. It follows that $V(f)\subset V$. Conversely, let $x\in V$. For each $j\ge1$, we take $C_j\in\mathcal{C}_{\rm ter}(f)$ and $m_j\ge 0$ such that
\[
x\in f^{-m_j}(U_{j,C_j}).
\]
Then, there are a sequence $1\le j_1<j_2<\cdots$ and $C\in\mathcal{C}(f)$ such that
\[
\lim_{l\to\infty}C_{j_l}=C
\]
in $\mathcal{C}(f)$. Note that for every $\epsilon>0$, we have
\[
C_{j_l}\subset U_\epsilon(C)
\]
for all sufficiently large $l\ge1$. For every $l\ge1$, by
\[
f^{m_{j_l}}(x)\in U_{{j_l},C_{j_l}},
\]
we have
\[
\omega^\ast(x,f)=\omega^\ast(f^{m_{j_l}}(x),f)\subset U_{\epsilon_{j_l}}(C_{j_l}).
\]
By
\[
\lim_{l\to\infty}\epsilon_{j_l}=0,
\]
we obtain
\[
\omega^\ast(x,f)\subset U_{2\epsilon}(C)
\]
for all $\epsilon>0$, thus $\omega^\ast(x,f)\subset C$. From Lemma 1.1, it follows that $C=C(x,f)\in\mathcal{C}_{\rm ter}(f)$, implying $x\in V(f)$. Since $x\in V$ is arbitrary, we conclude that $V\subset V(f)$, proving the claim.

Next, we show that $W(f)$ is a dense $G_\delta$-subset of $X$. Since $V(f)$ is a dense $G_\delta$-subset of $X$, it suffices to show that $W(f)$ is a dense $G_\delta$-subset of $V(f)$. Letting
\[
W=\bigcap_{j\ge1}\bigcap_{m\ge0}\{x\in V(f)\colon C(x,f)\subset U_{\frac{1}{j}}(\{f^i(x)\colon i\ge m\})\},
\]
we have $W=W(f)$. Let 
\[
W_{j,m}=\{x\in V(f)\colon C(x,f)\subset U_{\frac{1}{j}}(\{f^i(x)\colon i\ge m\})\}
\]
for all $j\ge1$ and $m\ge0$. Given any $x\in W_{j,m}$, $j\ge1$, $m\ge0$, by compactness of $C(x,f)$, there are $0<r<\frac{1}{j}$ and $n\ge m$ such that
\[
C(x,f)\subset U_r(\{f^i(x)\colon m\le i\le n\}).
\]
We take $\epsilon>0$ with $r+2\epsilon<\frac{1}{j}$. Since $x\in V(f)$ and so $C(x,f)\in\mathcal{C}_{\rm ter}(f)$, there is $a>0$ such that $d(x,y)<a$ implies
\[
C(y,f)\subset U_\epsilon(C(x,f))
\]
for all $y\in X$. By continuity of $f$, we have $b>0$ such that $d(x,y)<b$ implies
\[
\{f^i(x)\colon m\le i\le n\}\subset U_\epsilon(\{f^i(y)\colon m\le i\le n\})
\]
for all $y\in X$. It follows that $d(x,y)<\min\{a,b\}$ implies
\[
C(y,f)\subset U_{r+2\epsilon}(\{f^i(y)\colon m\le i\le n\})\subset U_{\frac{1}{j}}(\{f^i(y)\colon m\le i\le n\})\subset U_{\frac{1}{j}}(\{f^i(y)\colon i\ge m\})
\]
for all $y\in X$. Since $x\in W_{j,m}$ is arbitrary, $W_{j,m}$ is an open subset of $V(f)$. Since $j\ge1$ and $m\ge0$ are arbitrary, we conclude that $W$ is a  $G_\delta$-subset of $V(f)$. It remains to prove that $W$ is a dense subset of $V(f)$. Let $j\ge1$ and $m\ge0$. Given any $x\in V(f)$ and $\epsilon>0$, since $C(x,f)\in\mathcal{C}_{\rm ter}(f)$, there is $0<a<\epsilon/2$ such that $d(x,y)<2a$ implies
\[
C(y,f)\subset U_{\frac{1}{3j}}(C(x,f))
\]
for all $y\in X$. Since $f$ has the shadowing property, we see that
\[
C(x,f)\subset U_{\frac{1}{3j}}(\{f^i(p)\colon i\ge m\})
\]
for some $p\in X$ with $d(x,p)<a$. By compactness of $C(x,f)$, we obtain
\[
C(x,f)\subset U_{\frac{1}{3j}}(\{f^i(p)\colon m\le i\le n\})
\]
for some $n\ge m$. By continuity of $f$, we have $b>0$ such that $d(p,q)<b$ implies
\[
\{f^i(p)\colon m\le i\le n\}\subset U_{\frac{1}{3j}}(\{f^i(q)\colon m\le i\le n\})
\]
for all $q\in X$. Since $V(f)$ is a dense subset of $X$, we have $d(p,q)<\min\{a,b\}$ for some $q\in V(f)$. Note that
\[
d(x,q)\le d(x,p)+d(p,q)<2a<\epsilon.
\]
It follows that
\[
C(q,f)\subset U_{\frac{1}{3j}}(C(x,f))\subset U_{\frac{1}{j}}(\{f^i(q)\colon m\le i\le n\})\subset U_{\frac{1}{j}}(\{f^i(q)\colon i\ge m\}),
\]
implying $q\in W_{j,m}$. Since $x\in V(f)$ and $\epsilon>0$ are arbitrary, $W_{j,m}$ is an open dense subset of $V(f)$. Since $j\ge1$ and $m\ge0$ are arbitrary, we conclude that $W$ is a dense subset of $V(f)$, proving the claim. Thus, the theorem has been proved.
\end{proof}

We conclude with a remark on the proof.

\begin{rem}
\normalfont
\begin{itemize}
\item The proof shows that $V(f)$ and $W(f)$ are $G_\delta$-subsets of $X$ for every continuous map $f\colon X\to X$.
\item For any continuous map $f\colon X\to X$, we can show that if $f$ has the shadowing property, then
\[
V(f)=\{x\in X\colon\text{$C(\cdot,f)\colon X\to\mathcal{C}(f)$ is continuous at $x$}\}.
\]
By this, since $\mathcal{C}(f)$ is a compact metrizable space, we can show that $V(f)$ is a $G_\delta$-subset of $X$.
\item Let $f\colon X\to X$ be a continuous map and let $\xi=(x_i)_{i\ge0}$ be a sequence of points in $X$. For $\delta>0$, $\xi$ is called a {\em $\delta$-limit-pseudo orbit} of $f$ if $d(f(x_i),x_{i+1})\le\delta$ for all $i\ge0$, and
\[
\lim_{i\to\infty}d(f(x_i),x_{i+1})=0.
\]
For $\epsilon>0$, $\xi$ is said to be {\em $\epsilon$-limit shadowed} by $x\in X$ if $d(f^i(x),x_i)\leq \epsilon$ for all $i\ge 0$, and
\[
\lim_{i\to\infty}d(f^i(x),x_i)=0.
\]
We say that $f$ has the {\em s-limit shadowing property} if for any $\epsilon>0$, there is $\delta>0$ such that every $\delta$-limit-pseudo orbit of $f$ is $\epsilon$-limit shadowed by some point of $X$. When $f$ has the s-limit shadowing property, by Lemma 2.1, we can easily show that $W(f)$ is a dense subset of $X$.
\end{itemize}
\end{rem}


\begin{thebibliography}{99}

\bibitem{A1} E.\:Akin, The general topology of dynamical systems. Graduate Studies in Mathematics, 1. American Mathematical Society, Providence, R.I., 1993.

\bibitem{A2} E.\:Akin, On chain continuity. Discrete Contin. Dynam. Systems 2 (1996), 111--120.

\bibitem{AHK} E.\:Akin, M.\:Hurley, J.\:Kennedy, Dynamics of topologically generic homeomorphisms. Mem. Amer. Math. Soc. 164 (2003).

\bibitem{A} D.V.\:Anosov, Geodesic flows on closed Riemann manifolds with negative curvature. Proc. Steklov Inst. Math. 90 (1967), 235 p.

\bibitem{AH} N.\:Aoki, K.\:Hiraide, Topological theory of dynamical systems. Recent advances. North--Holland Mathematical Library, 52. North--Holland Publishing Co., 1994.

\bibitem{B} R.\:Bowen, Equilibrium states and the ergodic theory of Anosov diffeomorphisms. Lecture Notes in Mathematics, 470. Springer--Verlag, 1975.

\bibitem{BMR} W.\:Brian, J.\:Meddaugh, B.\:Raines, Shadowing is generic on dendrites. Discrete Contin. Dynam. Systems, Ser. S 12 (2019), 2211--2220.

\bibitem{BS} J.\:Buescu, I.\:Stewart, Liapunov stability and adding machines. Ergodic Theory Dynam. Systems 15 (1995), 271--290.

\bibitem{C} C.\:Conley, Isolated invariant sets and the Morse index. CBMS Regional Conference Series in Mathematics, 38. American Mathematical Society, Providence, R.I., 1978.

\bibitem{HH} M.W.\:Hirsch, M.\:Hurley, Connected components of attractors and other stable sets. Aequationes Math. 53 (1997), 308--323.

\bibitem{H} M.\:Hurley, Attractors: Persistence, and density of their basins. Trans. Amer. Math. Soc. 269 (1982), 247--271.

\bibitem{K} N.\:Kawaguchi, Generic and dense distributional chaos with shadowing. J. Difference Equ. Appl. 27 (2021), 1456--1481.

\bibitem{KMOK} P.\:Ko\'scielniak, M.\:Mazur, P.\:Oprocha, \L.\:Kubica, Shadowing is generic on various one-dimensional continua with a special geometric structure. J. Geom. Anal. 30 (2020), 1836--1864.

\bibitem{KOU} P.\:Krupski, K.\:Omiljanowski, K.\:Ungeheuer, Chain recurrent sets of generic mappings on compact spaces. Topology Appl. 202 (2016), 251--268.

\bibitem{MO} M.\:Mazur, P.\:Oprocha, S-limit shadowing is $C^0$-dense. J. Math. Anal. Appl. 408 (2013), 465--475.

\bibitem{P1} S.Yu.\:Pilyugin, The space of dynamical systems with the $C^0$-topology. Lecture Notes in Mathematics, 1571. Springer--Verlag, 1994.

\bibitem{P2} S.Yu.\:Pilyugin, Shadowing in dynamical systems. Lecture Notes in Mathematics, 1706. Springer--Verlag, 1999.

\bibitem{PP} S.Yu.\:Pilyugin, O.B.\:Plamenevskaya, Shadowing is generic. Topology Appl. 97 (1999), 253--266.

\bibitem{W} S.\:Willard, General topology. Dover Publications, Inc., Mineola, NY, 2004.

\end{thebibliography}
\end{document}